\documentclass{amsart}
\usepackage{latexsym}
\usepackage{amsmath,amssymb,amsfonts,amsthm,eufrak}
\usepackage{mathrsfs}
\usepackage{upref}
\usepackage{amscd}
\usepackage{eucal}
\usepackage[all]{xy}
\usepackage[dvips]{graphicx}
\newcounter{qcounter}

\usepackage{graphicx}
\usepackage{amsmath}
\usepackage[english]{babel}
\usepackage{amsfonts}
\usepackage{delarray}
\usepackage{bbold}
\usepackage{amssymb}
\usepackage{mathrsfs}
\usepackage{mathtools}
\usepackage{hyperref}
\usepackage{yfonts}
\usepackage{fancybox}
\usepackage{times}
\usepackage{amscd}

\title[Jacobian Newton Polyhedra]{Jacobian Newton Polyhedra and equisingularity\\ \ \\ \small{(Lecture
 at the Kyoto Singularities Symposium,\\ R.I.M.S., April 10, 1978.)}}
\author{B. Teissier}
\date{}
\address{Institut Math\'{e}matique de Jussieu, UMR 7586 du CNRS, 175 Rue du Chevaleret, 75013 Paris, France}
\email{teissier@math.jussieu.fr}

\begin{document}

\newtheorem{definition}{Definition}[section]
\newenvironment{define}{\begin{definition}\em}{\end{definition}}
\newtheorem{remark}{Remark}[section]
\newtheorem{remarks}{Remarks}[section]
\newtheorem{ex}{Example}[section]
\newenvironment{example}{\begin{ex}\em}{\end{ex}}
\newtheorem{theorem}{Theorem}[section]
\newtheorem{corollary}{Corollary}[section]
\newtheorem{claim}{Claim}[theorem]
\newtheorem{proposition}{Proposition}[section]
\newtheorem{lemma}{Lemma}[section]

\newcommand{\la}{\lambda}
\newcommand{\Hc}{\mathscr{H}}
\newcommand{\ma}{\mathscr{m}}
\newcommand{\Oc}{\mathcal{O}}
\newcommand{\Tc}{\mathscr{T}}
\newcommand{\CP}{\mathbb{CP}^2}
\newcommand{\Pc}{\mathbf{P}}
\newcommand{\B}{\mathbf{B}}
\newcommand{\C}{\mathbf{C}}
\newcommand{\D}{\mathbf{D}}
\newcommand{\R}{\mathbf{R}}
\newcommand{\Q}{\mathbf{Q}}
\newcommand{\Z}{\mathbf{Z}}
\newcommand{\N}{\mathbf{N}}
\newcommand{\F}{\mathscr{F}}

\newlength{\szer}
\newcommand{\Teiss}[2]{%
\settowidth{\szer}{$\displaystyle\frac{#1}{#2}$}%
\setlength{\szer}{0.5\szer}%
\left\{\hspace{\szer}%
\raisebox{0.14ex}{\makebox[0pt]{$\displaystyle\frac{#1}{\phantom{#2}}$}}%
\raisebox{-0.14ex}{\makebox[0pt]{$\displaystyle\frac{\phantom{#1}}{#2}$}}%
\hspace{\szer}\right\}%
}

\begin{abstract} In order to obtain invariants of the geometry of germs of complex hypersurfaces ``up to equisingularity" we extend the map $\mu$ (``jacobian multiplicity") which to a hypersurface associates its Milnor number to a map $\nu_j $ (``jacobian Newton polygon") which takes values in a monoid which is, in a way, the simplest monoid after the integers: the monoid of Newton polygons.

We show below some geometric reasons to do this, and summarize some results which show that $\nu_j([X_0])$ unites together many invariants of hypersurfaces.

In an appendix, we sketch some of the structure of the monoid of Newton polyhedra, with emphasis on the polygon $(d=2)$ case, and show some of the connections with the theory of multiplicities, mixed multiplicities, and integral closure of ideals.

The results and notations of the appendix are used in the text.

\end{abstract}

\maketitle

\section{Introduction}

We are interested\footnote{\textit{ I am grateful to Lilia Alanis L\'opez and Carlos Guzman, from CIMAT in Guanajuato, Mexico, who transcribed my old handwritten notes into LaTeX.}} in finding invariants of the geometry ``up to equisingularity" (see below) of a complex-analytic hypersurface $f(z_1,...,z_n)=0$ in $\C^n$, in the neighborhood of a singular point which we assume to be $0\in \C^n$.

To obtain invariants, we combine two ideas:

\vspace{5pt}
 $(1)$ The first is to compare the foliation defined by $f$ (in a neighborhood of $0\in \C^n$) - that is, the foliation having as leaves the level hypersurfaces $f=t$, with the foliation defined by some ``known" function $g:(\C^n,0)\rightarrow (\C,0)$ for example a general linear function on $\C^n$: We can extend this idea to the use of ``known" foliations of codimension $\geq 1$, say defined by a regular sequence $(g_1,...,g_k)$ on $\C\{z_1,...,z_n\}$ that is, having leaves of dimension $n-k$ defined by $g_1=u_1,...,g_k=u_k$.

The most natural way of comparing two foliations is to study the space of points where the leaves are not in general position, that is, in our case, the zero-space of the differential form $df\wedge dg_1\wedge \cdots dg_k$. This space is defined by the ideal generated by the $(k+1)-$minors of the jacobian matrix of $(f,g_1,...,g_k)$ and is nothing but the critical subspace of the map: $$p:\C^n\rightarrow \C^{k+1}\textrm{ (coordinates $t,u_1,...,u_k$)}$$ defined by $t\circ p=f$, $u_i\circ p=g_i$ ($1\leq i\leq k$).

If we assume that $f=g_1=\cdots =g_k=0$ has an isolated singularity, this critical subspace $C$ is finite over $\C^{k+1}$ and therefore we can define its \emph{image space} $D\subset \C^{k+1}$ (by the Fitting ideal $F_0(p_\ast \Oc_C)$) which is a hypersurface in $\C^{k+1}$ (see \cite{T1}, \S 1). It is simpler to study $D$ than to study $C$, and to do it we introduce the second idea:

\vspace{5pt}

 $(2)$ The second idea is to study a hypersurface by means of its \emph{Newton polyhedron}. Recall that the Newton polyhedron of a power series, say $$H(v_1,...,v_N)=\sum_{A\in\N^N}c_Av^A$$ is defined to be the (boundary of the) convex hull in $\N^N$ of the set $\bigcup_{c_A\neq 0}(A+\N^N)$.

The Newton polyhedron is not at all invariant by change of the coordinates $v_1,...,v_N$ and one must always specify which coordinates one chooses to compute it.

The procedure we propose is to take as invariants of $f$ the collection of Newton polyhedra of the equations defining the discriminants $D$ of maps $p$ as in $(1)$, with respect to the specified coordinates $(t,u_1,...,u_k)$ and for \emph{well chosen} sequences $(g_1,...,g_k)$.

 Here are three basic examples:

\begin{list}{\alph{qcounter})~}{\usecounter{qcounter}}

\item Case $k=n$, $g_i=z_i$ $(1\leq i\leq n)$: then our map $p: \C^n\rightarrow \C^{n+1}$ and its discriminant is nothing but its image, which is the hypersurface in $\C^{n+1}$ defined by: $t-f(u_1,...,u_n)=0$. In this case, our $D$ is the graph of $f$, and its Newton polyhedron is obtained directly from the Newton polyhedron of $f$ in the coordinates $z_1, ...,z_n$ in a way which the reader can see immediately . So this case is only the direct application of $(2)$ to $f$ without going through $(1)$.

\item Case $k=1,\ g_1=z_1$, a general linear function on $\C^n$ taken as coordinate; here ``general" means (\textit{cf.} \cite{C.E.W.}, Chap II) that the direction of the hyperplane $z_1=0$ is not a limit direction of tangent hyperplanes to $f=0$ at non-singular points. For simplicity, we assume that $f(z_1,...,z_n)=0$ has an isolated singularity. In any case what we have to study now is the plane curve $D\subset \C^2$, discriminant of the map $p:\C^n\rightarrow \C^2$ defined by $t=f$, $u_1=z_1$. \emph{From now on, I use freely the notations and results of the Appendix on Newton polygons at the end of this note.}

The Newton polygon of our plane curve $D$ in the coordinates $(t,u_1)$ is given by 
$$N(D)=\sum_{q=1}^l\Teiss{e_q+m_q}{m_q},$$
\noindent where $l$  is the number of irreducible components of the curve in $\C^n$ defined by $(\frac{\partial f}{\partial z_2},...,\frac{\partial f}{\partial z_n})$, $m_q$ is the multiplicity at 0 of the \textit{q}-th component of that curve, and $e_q+m_q$ is the intersection multiplicity at 0 of the same \textit{q}-th component with the hypersurface $f=0$. Note that by elementary properties of the intersection multiplicity, this last intersection multiplicity is $\geq m_q$, hence $e_q\geq 0$. In fact it can be shown that $e_q\geq m_q$, (see [CEW] chap II) because $e_q$ is the intersection multiplicity of the $q$-th component of the curve in $\C^n$ defined by the ideal $(\frac{\partial f}{\partial z_2},...,\frac{\partial f}{\partial z_n})$ with the hypersurface $\frac{\partial f}{\partial z_1}=0$; the intersection multiplicity $(\Gamma ,H)_0$ of a curve $\Gamma$ and a hypersurface $H$ in a non-singular space always satisfies $(\Gamma ,H)_0\geq m_0(\Gamma)\cdot m_0(H)$, hence we have $e_q\geq m_q$.

\item Case $k=0$: we consider $f:(\C^n,0)\rightarrow (\C,0)$, and assume that its fibre $(X_0,0)=(f^{-1}(0),0)$ has an isolated singularity. In this case, the critical subspace is defined by the ideal  $(\frac{\partial f}{\partial z_1},...,\frac{\partial f}{\partial z_n})=j(f)$, and the discriminant is defined in $(\C,0)$ by $t^{\mu}\C\{t\}$ where $\mu = \dim_\C\C\{z_1,...,z_n\}/j(f)$ is the Milnor  number of $f$. Here, the ``Newton polygon" is given as the ``convex hull" of the set $\mu+\N$: it is just the point with abscissa $\mu$.


\end{list}

Looking at a), b) and c), we see that even in the special case where we take the $g_j$ to be linear functions, the set of Newton polyhedra we obtain contains in particular:

\begin{itemize}

\item The polyhedra in $\R^{n+1}=\R\times \R^n$ which are obtained as convex join of the point $(1,0,...,0)$ with a Newton polyhedron of $f$, in all possible coordinate systems.

\item Polygons in $\R^2$, about which we shall see more below.

\item The Milnor number of $f$.

\end{itemize}

Our idea is to consider the monoid of integers, where the invariant ``Milnor number of $f$" takes its values, as a special case of others monoids, namely the monoids of Newton polyhedra. One of the advantages of this viewpoint is that it gives a natural frame to extend the results for isolated singularities to the general case, but we will not discuss this here, and will restrict ourselves to the isolated singularity case.

\vspace{5pt}

\textbf{\S 1.} In this lecture, I only want to list some of the invariants which are obtained by considering the case b) above, i.e., $k=1$ and a general linear function. Obviously this is the simplest case after $k=0$, i.e., the invariant we obtain (if indeed it is an invariant...) is the simplest one after the Milnor number.

First of all we must make clear that the Newton polygon we study can be obtained, up to elementary transformations, in these different ways: first, of course there is the way given in b) above:

\textit{i)} $ \sum_{q=1}^l\Teiss{e_q+m_q}{m_q}$ is the Newton polygon of the discriminant of the map
$\C^n\rightarrow \C^2$ given by $t=f,\ u=z_1$ (in the coordinates $t,u$).
\vspace{5pt}

\textit{ii)} $\sum_{q=1}^l\Teiss{e_q+m_q}{m_q}$ is the Newton polygon of a \textsl{general vertical plane section of the discriminant of a miniversal unfolding of} $f$ (or of a versal deformation of $(X_0,0)=(f^{-1}(0),0)$, in the coordinates given by the natural decomposition of the miniversal unfolding space $\C^{\mu}=\C\times\C^{\mu -1}$. See \cite{T1}.

\vspace{5pt}

\textit{iii)} $\sum_{q=1}^l\Teiss{e_q}{m_q}$ is the Newton polygon of the image of the curve $\Gamma$ in $\C^n$ defined by the ideal  $(\frac{\partial f}{\partial z_2},...,\frac{\partial f}{\partial z_n})$ by the map $\varphi :\C^n\rightarrow\C^2$ given by $t=\frac{\partial f}{\partial z_1}$, $v=z_1$, the Newton polygon of the plane curve $\varphi_\ast (\Gamma)$ is taken in the coordinates $t$, $u$. [Here, the image $\varphi_\ast (\Gamma)$ is defined by the Fitting ideal $F_0(\varphi_\ast \Oc_\Gamma)$, see \cite{T1}. It is important that the images should be compatible with base change].

 We shall work with this last polygon $\sum \Teiss{e_q}{m_q}$, but clearly the datum of any of the three is equivalent to the datum of any other. We shall denote $\sum \Teiss{e_q}{m_q}$ by $\nu_j([X_0])$. The fact that it depends only upon $X_0$ follows from the results in the appendix which imply that (with the notations of the appendix):
 
 $$\nu_j([X_0])=\nu_{j(f)}(\mathbf{m}),$$
\noindent
where $j(f)=(\frac{\partial f}{\partial z_1},...,\frac{\partial f}{\partial z_n})$, $\mathbf{m}=(z_1,...,z_n)\subset \C\{z_1,...,z_n\}$.

The first basic fact is that $\nu_j([X_0])$ is in fact an invariant of equisingularity, for some notion of equisingularity.

\begin{definition} Let $(X_0,0)\subset (\C^n,0)$ be a germ of a complex-analytic hypersurface. For each integer $i$, $0\leq i\leq n$, consider the Grassmanian $G^{(i)}$ of \textit{i}-dimensional linear spaces through 0 in $\C^n$, and $\min_{H\in G^{(i)}}\mu (X_0\cap H, 0)$. It can be proved that the set of \textit{i}-planes for which this minimum is obtained is a Zariski open dense subset of $G^{(i)}$. The corresponding value is denoted by $\mu^{(i)}$ and called the Milnor number of a general $i$-plane section of $X_0$.
\end{definition}

 Note that it can be shown that $\mu^{(i)}$ is also the smallest possible Milnor number for the intersection of $X_0$ with a non-singular subspace of dimension $i$ through 0. In all this we agree that if the singularity is not isolated, its Milnor number is $+\infty$.

 In the same vein, one can easily define the \textbf{topological type of a general \textit{i}-plane section of }$X_0$ (as embedded in the $i$-plane).
 
\begin{definition} We call the \textbf{total topological type} of $(X_0,0)$ the datum of the topological type of all the general i-plane sections of $X_0$, ($0\leq i\leq n$).
\end{definition}

\begin{theorem}\label{main} {\rm (See \cite{T3})} Let $\pi :(X,0)\rightarrow (S,0)$ be a family of $(n-1)$-dimensional complex analytic hypersurfaces with isolated singularities, where $S$ is non-singular. Assume that each fibre $X_s=\pi^{-1}(s)$ contains a point $x(s)$ such that $(X_s,x(s))$ has the same total topological type as $(X_0,0)$. Then, for a sufficiently small representative of $\pi$, we have that $x(s)$ is the only singular point of $X_s$, that there exist a complex analytic section $\sigma :S\rightarrow X$ of $\pi$ with $\sigma (s)=x(s)$ and that: 

\begin{center}
\underline{$\nu_j([X_s])$ is independent of $s\in S$}
\end{center}
\end{theorem}

In short, if we have a family of hypersurfaces where the total topological type is constant, $\nu_j ([X_s])$ is constant.

It should be remarked here that any ``reasonable" notion of equisingularity should imply the condition above, so that our $\nu_j ([X_0])$ will be constant in an equisingular deformation, once we know what that  term should mean. We remark also that the condition ``total topological type constant" implies, and in fact is equivalent to the condition: the sequence $\mu^\ast (X_s)$ of the Milnor number of general $i$-plane section of $X_s$ at $x(s)$ is independent of $s$, i.e.,

 $$\mu^\ast (X_s)=(\mu^{(n)} (X_s),\mu^{(n-1)} (X_s),...,\mu^{(1)} (X_s),\mu^{(0)} (X_s))$$ is constant.

\vspace{5pt}

\textbf{\S 2.} Now here is a partial list of the geometric features of $X_0$ which one can read from $\nu_j([X_0])$ and which are therefore constant in a deformation where the total topological type is constant:

\textbf{1)} The length $\ell(\nu_j[X_0])$ is equal to $\mu^{(n)}(X_0)$, the Milnor number of $X_0$.

\textbf{2)} The height $h(\nu_j[X_0])=\mu^{(n-1)}(X_0)$, the Milnor number of a general hyperplane section of $X_0$.

\textbf{3)} Therefore $\ell(\nu_j[X_0])+h(\nu_j[X_0])=\mu^{(n)}(X_0)+\mu^{(n-1)}(X_0)$ is equal to the \textbf{diminution of class} which the presence of singularity isomorphic with $(X_0,0)$ imposes on a projective hypersurface.

 This means that if $X\subset\Pc^n$ is a projective hypersurface such that the closure in $\check{\Pc}^n$ of the set of points representing the tangent hyperplanes to $X$ at non-singular points is again a hypersurface $\check X\subset\check{\Pc}^n$ (the dual hypersurface of $X$) we have, setting $d={\rm degree\  of}\  X$, $\check d={\rm degree\  of}\  \check X$ and assuming that $X$ has only isolated singular points, the formula (\textit{cf.} \cite{T6}, Appendix 2):
 
 $$\check d=d(d-1)^{n-1}-\sum_{x\in Sing X}(\mu^{(n)}(X,x)+\mu^{(n-1)}(X,x))$$

\begin{remark} For a plane curve, it gives, remarking that $\mu^{(1)}(X_0)=m_0(X_0)-1${\rm :} $$\check d=d(d-1)-\sum_{x\in Sing X}(\mu (X,x)+m(X,x)-1)$$ where $\mu(X,x)$ is the Milnor number at $x$ and $m(X,x)$ is the multiplicity.
\end{remark}

So far, we have used from our Newton polygon only very trivial features, namely its height and length, and the gist of theorem \ref{main} is as follows: if we look at a family of hypersurfaces, and the corresponding family of image curves (as in \textit{iii)} of \S 1), we know that the Newton polygons of the fibers of that family of curves all have the same length and height since these are $\mu^{(n)}$ and $\mu^{(n-1)}$ respectively. But this by no means implies in general that the Newton polygon is constant. For example consider the family of plane curves (depending on $\la$) $$t^a-u^b+\la t^p u^q=0\ \  \textrm{for given}\  a, b, {\rm say }\geq 3,\ {\rm where}\ p>0,q>0,\ \frac{p}{a}+\frac{q}{b}<1.$$

For all values of $\la$ the Newton polygon has height $a$ and length $b$, but it is not constant.

We see that our family of curves is of a rather special type, since as soon as height and length are constant the whole polygon is constant. Probably this is linked with the fact that the normalization of the discriminant of a versal unfolding is non-singular.

Anyway, now come some results which do use slopes of edges of the jacobian Newton polygon, albeit mostly that of the last edge:
 
\textbf{4)} Given $f(z_1,...,z_n)$ with isolated singularity at 0 and denoting by $m$ the maximal ideal of $\C\{z_1,...,z_n\}$, consider the invariant $\sum\Teiss{e_q}{m_q}$ associated with $f$ as in \S 1, \textit{iii)}. Then we have:
 
 For an integer N, the following properties are equivalent:
 \begin{enumerate}
 \item $N>\sup_q (\frac{e_q}{m_q})$
 
 \item Any function $g\in \C\{z_1,...,z_n\}$ such that $g-f\in m^{N+1}$ has the same topological type as $f$ (as a germ of mapping $(\C^n,0)\rightarrow (\C,0)$).
  \end{enumerate}
  
\textbf{5)} The smallest possible exponents in the \L ojasiewicz inequalities

 $$|\textrm{grad} f(z)|\geq C_1|f(z)|^{\theta_1}, \textrm{ }|\textrm{grad} f(z)|\geq C_2|z|^{\theta_2}\textrm{ (as } z\rightarrow 0)$$
 
 are given by: $\theta_1=\sup_q(\frac{e_q}{e_q+m_q})$, $\theta_2=\sup_q(\frac{e_q}{m_q})$.
 
\textbf{6)} In the note \cite{T2}, I introduced the following invariant of $f$: $\delta(f)$, or $\delta (X_0,0)$, is the maximum number of singular points which can appear in \underline{the same fibre} of an arbitrary small perturbation of $f$ (resp. deformation of $X_0$). Mr. I.N. Iomdin showed that.

$$\sup_q(\frac{e_q}{m_q})< 2\delta \leq \mu^{(n)}+\mu^{(n-1)}$$

In fact, he showed that as long as $k<\sup_q(\frac{e_q}{m_q})$, one can find arbitrarily close to $X_0$, a singularity of type $A_k$ i.e., isomorphic to: $z_1^{k+1}+z_2^2+\cdots +z_n^2=0$.

In modern notation:
$$k<\sup (\frac{e_q}{m_q})\Longrightarrow " A_k\longleftarrow X_0 \textrm{ (arrow=generalization).}$$

(hence $2\delta (X_0)\geq 2\delta (A_k)$ which is $k$ if $k$ is even, $k-1$ if $k$ is odd. However, this bound is not the best possible, i.e., one cannot reverse the implication above: For example the singularity $E_6$ ($z_1^2+z_2^3+z_3^4=0$) has all $\frac{e_q}{m_q}$ equal to 3 by the results in \cite{T3} but we have $E_6\rightarrow A_5$.

Since the $\frac{e_q}{m_q}$ of $A_5$ are all equal to 5, this also shows that \underline{there is no upper or lower} \underline{semi-continuity of} $\sup_q(\frac{e_q}{m_q})$ in a family of hypersurfaces, in general. (Note that $\nu_j(X_t)$ is upper-semi-continuous in the sense of polygons ``being above".) However , I do not know a counter-example to the possibility that: 

``In a family where $\mu^{(n)}$ is constant, $\sup_q (\frac{e_q}{m_q})$ is also constant".

(All the unimodular families are, up to suspension, families of curves with $\mu^{(2)}$ constant, hence $\mu^\ast$ constant, hence constant jacobian Newton polygon.)

\begin{remark} The invariant $\delta (X_0,0)$ is interesting for the following two reasons: first, in the case where $(X_0,0)$ is a plane curve, with ring $\Oc_0=\C\{z_1,z_2\}/(f(z_1,z_2))$ we have that: $\delta (X_0,0)=\dim_\C\overline{\Oc_0}/\Oc_0$ where $\overline{\Oc_0}$ is the integral closure of $\Oc_0$ in its total ring of fractions. Therefore, $\delta (X_0,0)$ is nothing but the ``diminution of genus" which the presence of a singularity isomorphic to $(X_0,0)$ would impose on a projective curve. Secondly, in any dimension, $\delta (X_0,0)$ gives a lower bound for the codimension of the Thom stratum of the origin in the base space of a versal unfolding. This result is one half of the ``Gibbs phase rule" (see \cite{T2}, \cite{T1}).
\end{remark}
\noindent
\textbf{Question:} Does the constancy of $\mu^{(n)}(X_t,0)$ imply the constancy of $\delta (X_t,0)$ for a family of hypersurfaces $(X_t,0)$?

For the case $n=2$, i.e., plane curves, the answer is yes, as was proved in \cite{T6} (Probably the result is still true for curves which are complete intersections, but here we deal only with hypersurfaces).

\begin{remark} From its definition and because the vertices of a Newton polygon have integral coordinates it is clear that $$\sup_q\frac{e_q}{m_q}\leq\mu^{(n)}(X_0).$$ so that if we like to use only the Milnor number for bounds, we may say that the $(\mu^{(n)}(X_0)+1)$-th Taylor polynomial of $f$ determines the topological type. However the result given above is much sharper; in fact, the equality $\sup\frac{e_q}{m_q}=\mu^{(n)}$ characterizes singularities of type $A_k$ (i.e., with $\mu^{(n-1)}=1$), as a glance at $\nu_j([X_0])$ will show. 
\end{remark}

\textbf{7)} Coming back to the study of the foliation of $\C^n$ by the level hypersurfaces $f=t$, we now describe some work of R\'emi Langevin.\par\noindent
Given a $C^\infty$-manifold of even dimension $V^{2n}\subset \R^M$ and a point $x\in V$, one will consider a normal vector $\mathbf{n}\in N_{V,x}$ and the orthogonal projection $\pi_\mathbf{n} : \R^M\rightarrow T\oplus \R.\mathbf{n}$, where $T_{V,x}\oplus \R.\mathbf{n}$ is the affine subspace of $\R^M$ through $x$ spanned by $T_{V,x}$ and $\mathbf{n}$ near $x$. The image $\pi_\mathbf{n}(V)$ is a hypersurface in this $2n+1$-dimensional affine space, having Gaussian total curvature $K_\mathbf{n}(x)$. We can take $\mathbf{n}$ of unit length, and the average of $K_\mathbf{n}(x)$ on the unit sphere of $N_{V,x}$ has a meaning. It is called the Lipschitz-Killing curvature of $V^{2n}$ and depends only upon the induced metric on $V^{2n}$. It is denoted by $K(x)$. (See \cite{L1} and \cite{L2}).

If we are given a Riemannian manifold $M$ foliated by even-dimensional submanifolds, we can define a map $K: M\rightarrow \R$ by: $K(m)=$ Lipschitz-Killing curvature at m of the leaf of the foliation going through $m$.

In the case of the foliations by $f=t$, the level manifolds are of even real-dimension, and all this is applicable, as follows:

Let $B_\epsilon =\{z\in\C^n\slash |z| \leq\epsilon\}$, let for $\zeta\in \R$, $\zeta>0$, $T_\zeta =\{z\in\C^n\slash |f(z)|<\zeta\}$, and $X_t=\{z\in \C^n \slash f(z)=t\}$.

\begin{proposition} {\rm (Langevin, see \cite{L1})} $$\lim_{\epsilon\rightarrow 0}\lim_{t\rightarrow 0}\int_{X_t\cap B_\epsilon}|K|dv=B_n\cdot(\mu^{(n)}(X_0)+\mu^{(n-1)}(X_0))$$ where $B_n$ is a constant  depending only on $n$. ($dv$= volume element)
\end{proposition}\noindent
(The double limit must be understood to mean that $\vert t\vert\rightarrow 0$ "much" faster than $\vert\epsilon\vert$).

\begin{proposition} {\rm (Langevin, see \cite{L2})}. For a fixed $\epsilon$, small enough, the integral $\int_{T_\zeta\cap B_\epsilon}|K|dw$ has a Puiseux expansion as a function of $\zeta$, and the first term is given by:

$$\int_{T_\zeta\cap B_\epsilon}|K|dw=\gamma_\epsilon \cdot \zeta^{2\frac{m_{q_0}}{e_{q_0}+m_{q_0}}}+\cdots\ ,$$
\noindent
where $\frac{e_{q_0}}{m_{q_0}}=\sup_q\frac{e_q}{m_q}$, and $\gamma_\epsilon\in\R$ is a positive real number ($dw$: volume element on $T_\zeta$).
\end{proposition}

The proofs of both results use the fact that the intersection multiplicity of the polar curve with respect to a general hyperplane direction, which is our curve $\Gamma$ defined by $(\frac{\partial f}{\partial z_2},...,\frac{\partial f}{\partial z_n})$, if $z_1=0$ is a general hyperplane, is equal to $\mu^{(n)}+\mu^{(n-1)}$, and that the $e_q$, $m_q$ appear in the parametrization of the branches of $\Gamma$. (See \cite{C.E.W.}, ch. II, \cite{I.E.P.}, ch. III)

\textbf{8)} Next, the filtration on the relative homology $H_{n-1}(X_t, X_t\cap H,\Z)$ [where $X_t:f=t$, $H:z_1=0$, where $H$ is general, and we look at everything in a small ball $B_\epsilon$] which is described in \cite{I.E.P.}, ch. III suggests that there could be some connection between the exponent  of oscilation associated to the critical point of $f$ at the origin, and the sequence of jacobian Newton polygons of $f$ and its restrictions to general $i$-planes through the origin $(1\leq i\leq n)$. Indeed, since it was proved that all the general hyperplane sections of $X_0$ are (c)-cosecant (see \cite{T3}) the jacobian Newton polygon of a general hyperplane sections is well defined, and then we can go on to lower dimensions.

Let $\gamma (t)$ be a horizontal family of homology classes of dimension $n-1$ in $X_t$ $(t\in\D_\eta$, $\eta$ sufficiently small), and let $\omega\in \Gamma (\C^n,\Omega^{n-1}_{\C^n/\D_\eta})$ [where $\C^n$ stands for $f^{-1}(\D_\eta)\cap\B_\epsilon$]. Consider the integral
$$I(t)=\int_{\gamma(t)}\omega$$ (\textit{cf.} the paper \cite{M} of Malgrange)

Then for small $t$ there is an expansion 
$$\int_{\gamma(t)}\omega=\sum_{ \alpha\in\Q,
 0\leq q\leq n-1}c_{\alpha, q}(\gamma,\omega)t^\alpha(\log t)^q $$ and the lower bound of the set of $\alpha$'s such that there exist $\gamma$, $\omega$, $q$ with $c_{\alpha,q}\neq 0$ is called the exponent of the Gauss-Manin connection of $f$, \textit{\`a la} Arnol'd, and denoted by $\sigma (f)$.
 
 Now let us denote by $\theta^{(i)}$ the rational number $\sup_q(\frac{e_q}{m_q})^{(i)}$, where $\nu_j\left[X_0\cap H^{(i)}\right]= \sum \Teiss{e_q}{m_q}^{(i)}$ is the jacobian Newton polygon of the intersection of $X_0$ with a general $i$-plane $H^{(i)}$ through $0$ in $\C^n$.\par\noindent
  \textbf{Question:} Does one have the inequality $$\sigma (f)\geq \sum_{i=1}^n\frac{1}{\theta^{(i)}+1}\textrm{ ?}$$ The underlying idea is that there should be a rather simple connection between $\sigma (f)$ and $\sigma (f|H)$ where $H$ is a general hyperplane through $0$. A more optimistic version of the same question is the following:
\par\noindent
 \textbf{Question:} Is there an expression $$\sigma(f)-\sigma(f|H)=\phi (\nu_j([X_0]))$$ where $\phi$ is a function from the set of Newton polygons to the rationals?
 
 In view of the fact that the behavior of $\sigma(f)$ under the Thom-Sebastiani operation is known, i.e., if $f\oplus f'$ is $f(z_1,...,z_n)+f'(z_{1}',...,z_{n'}')$ then $\sigma (f\oplus f')=\sigma(f)+\sigma(f')$, we are led to the following
 \par\noindent
 \textbf{Question:} Can one compute the sequence $\nu_j^\ast([X_0]\perp[X'_0])$ of the jacobian Newton polygons $\nu_j(([X_0]\perp[X'_0])\cap H^{(i)})$ of the sections by general $i$-dimensional planes of the hypersurface defined by $f\oplus f'=0$ $(1\le i\leq n+n')$ from the corresponding sequences $\nu_j^\ast([X_0])$ and $\nu_j^\ast([X'_0])$?
 
\textbf{9)} Finally we describe the only result so far which uses the totality of the jacobian Newton polygon, and which shows that at least in special case, the jacobian Newton polygon contains a lot of information.

Let $(X_0,0)\subset (\C^2,0)$ be a germ of an irreducible plane curve, given parametrically by $z_1=z_1(t)$, $z_2=z_2(t)$. Then $\Oc_{X_0,0}=\C\{z_1(t),z_2(t)\}\hookrightarrow\C\{t\}$ and it is known (\textit{cf.} \cite{Z}) that the topological type (or equisingularity type) is completely determined by, and determines, the semi-group $\Gamma$ of the orders in $t$ of the elements of $\Oc_{X_0,0}$.

Let $\Gamma=\langle \bar{\beta_0},\bar{\beta_1},...,\bar{\beta_g}\rangle$ be a minimal system of generators for this semi-group. Set $l_i=(\bar{\beta_0},...,\bar{\beta_i})$ (greatest common denominator and define $n_i$ by $l_{i-1}=n_i\cdot l_i$. Since $l_g=1$ we have $l_0=\bar{\beta_0}=n_1\cdots n_g$.

\begin{theorem} {\rm (M. Merle \cite{Me})}
$$\nu_j([X_0])=\sum_{q=1}^g\Teiss{(n_q-1)\bar{\beta_q}-n_1\cdots n_{q-1}(n_q-1)}{n_1\cdots n_{q-1}(n_q-1)},$$ \end{theorem}
\noindent
and therefore $\nu_j([X_0])$ is completely determined by $\Gamma$ and determines it: in this case, $\nu_j([X_0])$ is a \underline{complete invariant} of the topological type of $(X_0,0)$.\par\medskip\noindent
\textbf{Jacobian Newton Polygons and quasi-homogeneous singularities:}

One should be aware of the fact that the jacobian Newton polygon is \underline{not} determined by the weights in a quasi-homogeneous hypersurface with isolated singularity: the following example was found by Brian\c con and Speder:

Let $(X_t)$ $(t\in \C)$ be the family of surfaces defined by the equation: $z_2^3+tz_1^\alpha z_2+z_1^\beta z_3+z^{3\alpha}_3$ where $\alpha$, $\beta$ are given integers such that $\alpha\geq 3$ and $3\alpha =2\beta +1$. 

Then $$\nu_j([X_0])=\Teiss{2\beta}{2}+\Teiss{2\beta(2\beta -2)}{2\beta -2},$$

$$\nu_j([X_t])=\Teiss{2\beta (2\beta -1)}{2\beta -1} \textrm{ for }t\neq 0.$$ 

\begin{remark} In this example the usual Newton Polyhedron in the coordinates $z_1$, $z_2$, $z_3$ is not the same fot $t=0$ and $t\neq 0$, however, the plane supporting the only compact face remains the same of course.
\end{remark} 

 An optimistic question is the following:
 \par\noindent
 \textbf{Question:} Given a hypersurface consider the totality of its Newton polyhedra with respect to all coordinate systems. From each remember only the support hyperplanes of the compact faces. Is the set of possible configurations of these support hyperplanes determined by the topology (or the ``equisingularity class") of the hypersurface?
  
Check it for the case $n=2$.
\vfill\eject
\section{Appendix on Newton Polyhedra}

Let K be a field, and $f=\sum_{A\in \mathbf{N}^k}f_A z^A \in K \left[\left[x_1,\ldots,x_k\right]\right]$. Define $\mathrm{Supp}(f)=\{A\in\mathbf{N}^k | f_A\neq 0\},\:$ consider the set $\Delta(f)=\bigcup _{A\in \mathrm{Supp}(f)}(A + \mathbf{N}^k)\:$ and $N(f)\:$ its convex hull in $\mathbf{R}^k,\:$ which we call the Newton Polyhedron (in fact the boundary of this convex region is also called the Newton Polyhedron) of the series f.
Let $N_c(f)\:$ denote the union of the compact faces of $N(f).$\begin{remark} Assume that $f\in K\left[x_1,\ldots,x_k\right];\:$ the definition above is well adapted to the study of $f\:$ near the origin. Historically people also considered the convex hull of $\mathrm{Supp}(f),\:$ which when $f\:$ is a polynomial is a compact convex body in $\mathbf{R}^k,\:$ and used it to obtain information on the singularities of the hypersurface $\{f=0\}\:$ not only at the origin, but at the other singular points, including those at infinity.
\end{remark}
Let us define the sum of two Newton polyhedra $N_1\:$ and $N_2\:$ as follows: $N_1+N_2\:$ is the (boundary of the) convex hull in $\mathbf{R}^k\:$ of the set of points of the form $p_1+p_2$ where $p_i$ lies in the convex region of $\R^k$ bounded by $N_i\: (i=1,2).\:$ Then the sum is clearly commutative and associative.\\
Exercise 1: Show that if $f_1,\:f_2\: \in K\left[\left[x_1,\ldots,x_k\right]\right]\:$ then $N(f_1f_2)=N(f_2)+N(f_2).$\par\noindent
Exercise 2: Define a Newton Polyhedron to be \underline{elementary} if it has at most one compact face of dimension $k-1,\:$ and its non compact faces of dimension $k-1\:$ lie in coordinate hyperplanes, i.e., it is the standard simplex up to affinity. Say that a Newton polyhedron has \underline{finite volume} if the volume of the complement in the positive quadrant of the convex region it bounds is finite. Give an example of a polyhedron of finite volume which is \emph{not} a sum of elementary polyhedra.\par\noindent
Our viewpoint is to consider that Newton polyhedra, with the addition described above, constitute a very natural generalization of the monoid of integers.\par
 \subsection{ Study of the case $k=2$ i.e., Newton Polygons.}
Historically Newton first used the Newton Polygon to describe a successive approximation procedure for the computation of the roots $y=y(x)\:$ (fractional power series) of an algebraic equation $f(x,y)=0\:$ near $0,\:$ assumed to be an isolated singular point of $f\in K\left[\left[x,y\right]\right],\:K=\mathbf{R}\:{\rm or}\:\mathbf{C}.$\par\noindent
Before describing Newton's result, we need some notation: Any Newton polygon enclosing a finite area determines a length $\ell(N)\:$ which is the length of its projection on the horizontal axis (by convention) and a height $h(N),\:$ which is the length of its projection on the vertical axis:\par\medskip\noindent
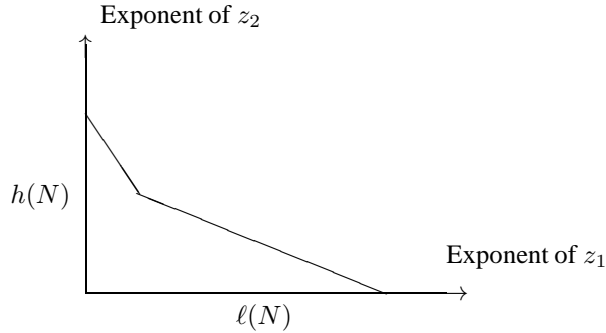
\vspace{.4in}
\begin{picture}(0,1)\setlength{\unitlength}{.4cm}
\put(7.5,-6){$h(N)$}
\put(10,-9){\line(0,1){8.29}}
\put(9.78,-1){$\uparrow$}
\put(10.5,0){Exponent of $z_2$}
\put(10,-9){\line(1,0){12}}
\put(22,-8){Exponent of $z_1$}
\put(10,-3){\line(2,-3){1.8}}
\put(20,-9){\line(-5,2){8.3}}
\put(21.85,-9.22){$\rightarrow$}
\put(15,-10){$\ell(N)$}
\end{picture}
\vspace{1.6in}
\vfill\eject
We need a notation for an elementary polygon as follows: an elementary polygon is entirely described by its length and height, and we write it $\Teiss{\ell(N)}{h(N)}$:

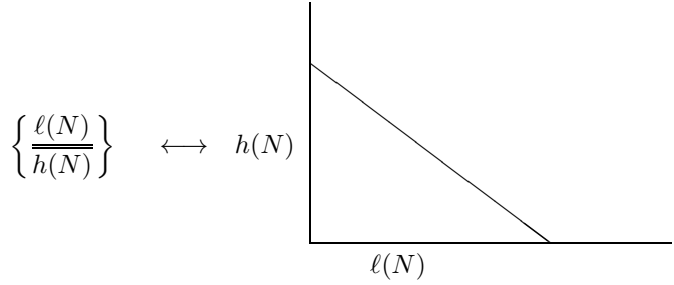
\begin{picture} (1,1)\setlength{\unitlength}{.4cm}
\put(4,-7){$\Teiss{\ell(N)}{h(N)}$}
\put(9,-7){$\longleftrightarrow$}
\put(11.5,-7){$h(N)$}
\put(14,-10){\line(0,1){8}}
\put(14,-10){\line(1,0){12}}
\put(14,-4){\line(4,-3){8}}
\put(16,-11){$\ell(N)$}
\end{picture}
\vspace{2in}

Let $N$ be a Newton polygon enclosing a finite area: it can be written (non-uniquely) as a sum of elementary polygons:

$$N=\sum_{i=1}^{\textit{l}}\Teiss{\ell(N_i)}{\textit{h}(N_i)}.$$

However, if we require that $i\neq j \Longrightarrow \frac{\ell(N_i)}{h(N_i)}\neq \frac{\ell(N_j)}{h(N_j)}\:$ the decomposition becomes unique, and we call it the canonical decomposition of $N$.\\

\subsection{ First use of the Newton Polygon: The Newton-Puiseux Theorem.}

Consider $f\in \mathbf{C}\{z_1,z_2\}\:$ (or $\:\mathbf{C}\left[\left[z_1,z_2\right]\right]$) such that $f(0,z_2)\:$ is not the function zero. By the Weierstrass preparation theorem, up to multiplication by a unit, which changes neither the germ at $0\:$ of the zero set $\{f=0\}\:$ no the Newton polygon, we may assume that $f \in \mathbf{C}\{z_1\}\left[z_2\right]$.\\

Since $\mathbf{C}\{z_1\}\:$ is an henselian ring, if we denote by $K\:$ its field of fractions and by $v\:$ the $(z_1)$-adic valuation, $v\:$ has a unique extension to every finite extension $L\:$ of $K$, [and this extension takes its values in $\frac{1}{d}\mathbf{Z}\:$ where $d=\left[L:K\right],\:$ by a theorem of Puiseux, but we do not need this, only the fact that it takes values in a submonoid of $\mathbf{Q}\:$ isomorphic with $\mathbf{Z}\:$].\\

Let $m_{\rho}\:$ be the number of roots of $f=0\:$ in a splitting extension $L\:$ for $f,\:$ which have a given valuation $\rho$.

\begin{theorem}[Newton-Puiseux]

					\begin{enumerate}
					\item[]
					\item[]
					\item[I.]							    							$N(f)=\sum_{\rho\in\mathbf{Q}} 			\Teiss{m_{\rho}\rho}{m_{\rho}}.$
			
			\item[]
			\item[II.]If $f\:$ is irreducible in $\mathbf{C}\{z_1,z_2\}\:$ or $\mathbf{C}\{z_1\}\left[z_2\right],\: N(f)\:$ is elementary, and the canonical decomposition of $N(f)\:$ into elementary polygons \underline{can be realized} by a decomposition $f=f_1\ldots f_k$.\\
				\end{enumerate}
	\end{theorem}

\begin{corollary}
Let $U\subset \mathbf{C}\{z_1\}\left[z_2\right]\:$ be the submonoid (for the multiplication) consisting of unitary polynomial in $z_2,\:$ not divisible by $z_2$ and let $NU\:$ the set of Newton polygons which meet both axis, a submonoid of the additive monoid of all the Newton polygons (see the second exercise after definition \ref{special} below) and let $\sim$ be the equivalence relation on $U\:$ defined by: $p_1\:\sim\:p_2\:$ if $p_1\:$ and $p_2\:$ have the same number $m_{\rho}\:$ of roots with a given valuation $\rho\:$ in any finite Galois extension $L\:$ of $K$. Then $f\longrightarrow N(f)\:$ induces an isomorphism of monoids $U/\sim\: \longrightarrow NU\:$ (remark that $\sim\:$ is compatible with the multiplication on $U$). 
\end{corollary}

\subsection{ Second use of the Newton Polygon: The Newton Polygon of two ideals.}

\subsubsection{}
 Let $\Oc\:$ be a Cohen-Macaulay reduced analytic algebra, and let $\mathbf{n}_1\:$ and $\mathbf{n}_2\:$ be two ideals of $\Oc\:$ which are primary for the maximal ideal. Consider a representative $W\:$ of the germ of complex analytic space corresponding to $\Oc,\:$ and let 

$$\pi:\overline{W'}\longrightarrow W$$
\noindent
be the normalized blowing up of the ideal $\mathbf{n}_1.\:$ Since $\Oc\:$ is Cohen-Macaulay this normalized blowing up is also the normalized blowing-up of an ideal $\mathbf{n}_1^{\left[d\right]}\subset\mathbf{n}_1\:$ where $d=\mathrm{dim} \Oc,\:$ and $\mathbf{n}_1^{\left[d\right]}\:$ is an ideal generated by $d\:$ \textquotedblleft sufficiently general\textquotedblright\ elements of $\mathbf{n}_1.\:$ Therefore $\mathbf{n}_1^{\left[d\right]}\:$ is generated by a regular sequence, say $(f_1,\ldots,f_d)\:$ and its blowing up $\pi_0: W'\:\longrightarrow W\:$  can be identified with the restriction of the first projection to the subspace $W'\subset W \times \mathbf{P}^{d-1}\:$ defined by the ideal $(f_iT_j-f_jT_i,\: 1\leq i < j\leq d)$. The exceptional divisor is set-theoretically $\mathbf{P}^{d-1}=\{0\}\times \mathbf{P}^{d-1}\subset W'\:$ (since $\mathbf{n}\subset\mathbf{m}\:$, the maximal ideal of $\Oc=\Oc_{W,0}$) and therefore, each irreducible component of the exceptional divisor $D\:$ of $\pi\:$ (defined by $\mathbf{n}_1.\Oc_{\overline{W'}}$) maps to a divisor on $W'\:$ which is set-theoretically $\mathbf{P}^{d-1}\:$:  the normalization map $n\colon\overline{W'}\longrightarrow W' $ induces a surjective map $D_{i,\mathrm{red}}\longrightarrow\mathbf{P}^{d-1}$. Let $\mathrm{deg} D_{i,\mathrm{red}}\:$ denote the degree of this map. Of course it can also be computed directly from $\mathbf{n}_1.\Oc_{\overline{W'}}$, as follows: setting $\mathcal{I}=\mathbf{n}_1.\Oc_{\overline{W'}},\:$ a $\pi$-ample invertible sheaf, we have 
$$\mathrm{deg} D_{i,\mathrm{red}}=\mathrm{deg}(\mathcal{I}/\mathcal{I}^2\otimes_{\Oc_D}\Oc_{D_i,\mathrm{red}}),$$
\noindent
where the degree of an ample invertible sheaf $\mathcal{L}\:$ on $X\:$ is given by:
$$\chi(\mathcal{L}^{\otimes\nu})=\frac{\mathrm{deg}\mathcal{L}}{(\mathrm{dim} X)!}\:\:\nu^{\mathrm{dim} X}+ O(\nu^{\mathrm{dim}X-1}).$$
\noindent
Now since $\overline{W'}\:$ is normal, and $D\:$ is a divisor, at a \textquotedblleft  general\textquotedblright\  point of each $D_{i,\mathrm{ red}}\:$ we can define the order of vanishing of an element of the form $f\circ \pi,\:f\in\Oc\:$ along $D_i$: at a general point $z\:$ of $D_i\:$ both $\overline{W'}\:$ and $D_{i,\mathrm{red}}\:$ are non-singular, and $D_i\:$ is defined by $v^{\nu}.\Oc_{\overline{W'},z}\:$ where $v\:$ is part of a local coordinate system. Define $v_{D_i}(f\circ\pi,z)\:$ as the highest power of $v\:$ dividing $f\circ\pi\:$ in $\Oc_{\overline{W'},z}$. \\

This number is obviously locally constant, hence depends only upon $D_i$, since $z\:$ lie in a Zariski open subset of $D_i$, which is connected since $D_i\:$ is irreducible.
Therefore we have associated with each component $D_i\:$ of $D\:$ an application 
$$	\begin{array}{lccl}
	 v_{D_i}:&\Oc&\longrightarrow&\mathbf{N}\\
	 				& f   &\longmapsto&v_{D_i}(f\circ\pi,z)\\
	\end{array}$$

where $z\:$ stays in some Zariski open subset of $D_i$.

We define, for an ideal $\mathbf{n}\:$ of $\Oc,\:$ $v_{D_i}(\mathbf{n})=\mathrm{min} _{h\in\mathbf{n}}\{v_{D_i}(h)\}$.

\begin{definition}
Let $\mathbf{n}_1\:$ and $\mathbf{n}_2\:$ be two primary ideals of $\Oc$. The Newton polygon of $\mathbf{n}_1\:$ and $\mathbf{n}_2\:$ is defined to be:
$$\nu_{\mathbf{n}_1}(\mathbf{n}_2)=\sum_{i=1}^{r}\mathrm{deg}D_{i,\mathrm{red}}\Teiss{v_{D_i}(\mathbf{n}_1)}{v_{D_i}(\mathbf{n}_2)}$$
where the exceptional divisor $D\:$ of the normalized blowing up $\pi:\overline{W'}\longrightarrow W\:$ of $\mathbf{n}_1\:$ is $D=\bigcup_{i=1}^r D_i$, with $D_i\:$ irreducible.
\end{definition}

Now let us consider an ideal $\mathbf{n}_1^{\left[d-1\right]}\:$ generated by $d-1\:$ \textquotedblleft  general\textquotedblright\   elements of $\mathbf{n}$. It defines a curve $\Gamma\subset W\:$ which is reduced. Let $\Gamma=\bigcup_{j=1}^{l}\Gamma_j\:$ where each $\Gamma_j\:$ is an irreducible germ of curve in $W$. By normalization we get a composed map $\overline{\Gamma_j}\stackrel{n_j}{\longrightarrow}\Gamma_j\subset W\:$ and therefore, since $\Oc_{\overline{\Gamma}_{j,0}}=\mathbf{C}\{t_j\}$, a new order function

$$	\begin{array}{lccl}
	 v_{\Gamma_j}:&\Oc&\longrightarrow&\mathbf{N}\\
	 				& f   &\longmapsto&\mathrm{order\:in\:}t_j\: \mathrm{of}\:(f\vert_{\Gamma_j}\circ n_j) \\
	\end{array}$$
for each irreducible component $\Gamma_j\:$ of $\Gamma$. Again we can define $v_{\Gamma_j}(\mathbf{n})\:$ for an ideal $\mathbf{n}\:$ in $\Oc$. Using the construction explained in \cite{T5} it is not difficult to prove:

\begin{proposition}
\begin{equation}\label{4.1.1}
\nu_{\mathbf{n}_1}(\mathbf{n}_2)=\sum_{j=1}^{l}\Teiss{v_{\Gamma_j}(\mathbf{n}_1)}{v_{\Gamma_j}(\mathbf{n}_2)}
\end{equation}
\end{proposition}

\begin{corollary}
We have the two equalities
$$\begin{array}{ccc}
				h(\nu_{\mathbf{n}_1}(\mathbf{n}_2))&=&e(\mathbf{n}_1^{\left[d-1\right]}+\mathbf{n}_2^{\left[1\right]})\\
				
\ell(\nu_{\mathbf{n}_1}(\mathbf{n}_2))&=&e(\mathbf{n}_1)=e(\mathbf{n}^{[d]}_1)\\					\end{array}
$$
where $e(\mathbf{n}_1^{\left[i\right]}+\mathbf{n}_2^{\left[d-1\right]})\:$ are the \textquotedblleft mixed multiplicities\textquotedblright \ introduced in {\rm (\cite{C.E.W.}, Chap. I)} and in particular $e(\mathbf{n}_1)=\sum_{i=1}^{r}\mathrm{deg}D_{i,\mathrm{red}}v_{D_i}(\mathbf{n}_1)$, which yields a \textquotedblleft projection formula\textquotedblright\  which I used, in the case $d=2$, in \cite{T4}.
\end{corollary}

Proof: $\sum_{i=1}^{l}v_{\Gamma_i}(\mathbf{n}_2)\:$ is equal to $\sum_{i=1}^{l}v_{\Gamma_i}(h)\:$ where $h\:$ is a \textquotedblleft general\textquotedblright\  element of $\mathbf{n}_2$, and this is equal to the intersection multiplicity of $\{h=0\}\:$ with $\Gamma$, which is $e(\mathbf{n}_1^{\left[d-1\right]}+\mathbf{n}_2^{\left[1\right]})\:$ where $\mathbf{n}_2^{\left[1\right]}=h.\Oc$.

\begin{corollary}
If $\overline{\mathbf{n}_1'}=\overline{\mathbf{n}_1}, \:\overline{\mathbf{n}_2'}=\overline{\mathbf{n}_2}$, then
$$\nu_{\mathbf{n}_1'}(\mathbf{n}_2')=\nu_{\mathbf{n}_1}(\mathbf{n}_2)$$
where the bar denotes the integral closure of ideals.
\end{corollary}

Proof: The ideals $\mathbf{n}_1,\mathbf{n}_1'\:$ determine the same normalized blowing up, hence
$v_{D_i}(\mathbf{n}_1')=v_{D_i}(\mathbf{n}_1)$. Also $v_{\Gamma_j}(\mathbf{n}_2')=v_{\Gamma_j}(\mathbf{n}_2)\:$ by the valuative criterion of integral dependence.\\

\subsubsection{ The Newton polygon of two ideals as a dynamic version of intersection multiplicities.}

We want to explain informally what type of information $\nu_{\mathbf{n}_1}(\mathbf{n}_2)\:$ contains compared to $e(\mathbf{n}_1)$. Suppose that we have reduced an intersection multiplicity problem to the computation of the multiplicity of an ideal $\mathbf{n}_1=(f_1,\ldots,f_d)\:$ generated by a regular sequence in a Cohen-Macaulay analytic algebra $\Oc$. Then by flatness, one may compute $e(\mathbf{n}_1)\:$ as the number of intersections of the hypersurface $f_d-v=0\:$ with the curve $\Gamma:f_1=\ldots=f_{d-1}=0$. If our $(f_1,\ldots,f_{d-1})\:$ are sufficiently general these will be transversal intersections of nonsingular varieties for $v\neq 0\:$ sufficiently small (i.e., $f_d-v=0\:$ meets $\Gamma\:$ at non-singular points of $\Gamma$, and transversally). Now on a component $\Gamma_j\:$ of $\Gamma$, we have an expansion for $f_d|_{\Gamma_j}\:$ say $f_d|_{\Gamma_j}=\gamma_j t_j^{e_j}+\dots\:(e_j\in\mathbf{N})\:$ where $t_j\:$ is the local coordinate on $\overline{\Gamma}_j$. Now let $\mathbf{m}\:$ be the maximal ideal of $\Oc$, and $m_j=v_{\Gamma_j}(\mathbf{m})$, so there is a $z_0\in\mathbf{m}\:$ such that $z_0|_{\Gamma_j}=\delta_jt_j^{m_j}+\dots$\\

This implies that we can express the value of $z_0\:$ at the intersection points of $\Gamma\:$ and $f_d-v=0\:$ by a Puiseux expansion \underline{as a function of $v$} on each $\Gamma_j$:

\begin{equation}\label{4.2.1}
z_0= \epsilon_j v^{\frac{m_j}{e_j}}+\dots
\end{equation}

In all there are $e_j\:$ intersection points on $\Gamma_j$, and $\sum_{j}e_j=e(\mathbf{n}_1)$. If we take $\Oc=\mathbf{C}\{z_0,\ldots,z_d\}$, the expansion (\ref{4.2.1}) above tells us that our Newton polygon
$$\nu_{\mathbf{n}_1}(\mathbf{m})=\sum_{j=1}^{r}\Teiss{e_j}{m_j}$$
describes precisely the collection of \textquotedblleft vanishing rates\textquotedblright \  of (a \textquotedblleft generic\textquotedblright\  coordinate of) the intersection points of $\Gamma\:$ with $f_d-v=0$, as functions of $v$, whereas $e(\mathbf{n}_1)\:$ is only the number of these points.\\

\subsubsection{The product formula.}

Let $\Oc\:$ and $\Oc'\:$ be two analytic algebras (in fact two noetherian local rings), and let $\mathbf{n}_1\subset\Oc$, $\mathbf{n}_1'\subset\Oc'\:$ be ideals which are primary for the respective maximal ideals. Then the multiplicity in the tensor product $\Oc\otimes_{\mathbf{C}}\Oc'\:$ of the ideal\break $N=\mathbf{n}_1\otimes 1+1\otimes \mathbf{n}_1'\:$ is given by:
$$e(N)=e(\mathbf{n}_1).e(\mathbf{n}_1').$$

From our viewpoint it is natural to ask whether one can have a similar product formula for the Newton polygons of two ideals, with a suitable definition of the product of two Newton polygons. The answer is yes:

\begin{definition}
Let $P=\Teiss{\ell(P)}{h(P)}\:$ and $Q=\Teiss{\ell(Q)}{h(Q)}\:$ be two elementary Newton polygons; define the product $P*Q\:$ by

$$P*Q=\Teiss{\ell(P) \ell(Q)}{\mathrm{min}(\ell(P)h(Q),\ell(Q)h(P))},$$
and for any two polygons $P,Q\:$ of finite volume, taking decompositions $P=\sum_i P_i$, $Q=\sum_j Q_j\:$ into elementary polygons, define 

$$P*Q=\sum_{i,j}P_i*Q_j.$$
(The result is independent of the choice of decompositions)
\end{definition}
\noindent
Exercise: Check that $*\:$ is commutative, associative and distributive (!)

\begin{proposition}
Let $\Oc$, $\Oc'\:$ be two Cohen-Macaulay analytic algebras, $\mathbf{n}_1$, $\mathbf{n}_2\:$ in $\Oc\:$ primary ideals, $\mathbf{n}_1'$, $\mathbf{n}_2'\:$ in $\Oc'\:$ also. Set $N_1=\mathbf{n}_1\otimes 1+1\otimes\mathbf{n}_1'$,  $N_2=\mathbf{n}_2\otimes 1+1\otimes\mathbf{n}_2'$, both ideals on $\Oc\otimes_{\mathbf{C}}\Oc'\:$ (note that we could just as well take the analytic tensor product $\Oc\stackrel{_\frown}{\otimes} _{\mathbf{C}}\Oc'\:$). Then we have:
$$\nu_{N_1}(N_2)=\nu_{\mathbf{n}_1}(\mathbf{n}_2)*\nu_{\mathbf{n}_1'}(\mathbf{n}_2').$$
\end{proposition}
The proof follows the construction given in \cite{T3}.

\subsubsection{Comments on the product of Newton polygons.}

One could say that the definition of the product of Newton polygons is just ad-hoc for the proposition above. However, it has some interesting features, of which I now show two:\\

First, we try to find an operation on polynomials which induces, bye taking Newton polygons, the product $*$. Let $V\:$ be a ring (here $\mathbf{C}\{z_1\}$) and consider the following operation on the ring $V\left[T\right]\:$ (here $T=z_2$).

$$P_1,\:P_2\in V\left[T\right]\rightsquigarrow \mathrm{Res}^{V}_{P_1,P_2}(T)=\mathrm{Res}_U(P_1(T+U),P_2(U)),$$
 \noindent
where $\mathrm{Res}_U\:$ means the \underline{resultant} with respect to the new indeterminate $U$. We assume now that $P_1\:$ and $P_2\:$ are unitary. We have that $\mathrm{Res}_{P_1,P_2}(T)\in V\left[T\right]\:$ is a polynomial of degree $\mathrm{deg} P_1. \mathrm{deg}P_2\:$ and $\mathrm{Res}^{V}_{P_1,P_2}(0)=\mathrm{Res}_{T}(P_1(T),P_2(T))$. If we go to a splitting extension for $P_1\:$ and $P_2\:$ of the fraction field of $V$, and let $\bar{V}\:$ be the integral closure of $V\:$ in it; $\mathrm{Res}^{\bar{V}}_{P_1,P_2}(T)=\Pi _{i,j}(T-(\alpha_i-\beta_j))\:$ where $\alpha_i,\beta_j\:$ are the roots of $P_1,\: P_2$.\\

\begin{proposition}
Let $N_1,\:N_2\:$ be two Newton polygons of finite volume. For \textquotedblleft almost all\textquotedblright\  pairs of elements $P_i\in\mathbf{C}\{z_1\}\left[z_2\right]\:$ such that $N(P_i)=N_i\:(i=1,2)\:$ we have that 
$$N(Res^{\mathbf{C}\{z_1\}}_{P_1,P_2}(z_2))=N_1*N_2.$$
\end{proposition}

For any finite edge $\gamma\:$ of the Newton Polygon of $P\in\mathbf{C}\{z_1\}\left[z_2\right]\:$ define $P_{\gamma}\:$ to be the sum of those terms in $P\:$ having an exponent which is on $\gamma$. Then \textquotedblleft almost all\textquotedblright\  means that for any pair $\gamma_1,\: \gamma_2\:$ of finite edges on $N_1,\:N_2$, the (quasi-homogeneous) polynomials $P_{1,\gamma_1}\:$ and $P_{2,\gamma_2}\:$ have no common zero outside $0$, (in a neighborhood of $0$). Then is clearly a Zariski-open condition on the coefficients of $P_{N_i}=\sum_{\gamma\subset N_i}P_{i,\gamma}$

The proof is straightforward from the expression of $\mathrm{Res}^{\bar{V}}_{P,l}(T)\:$ and the fact that the Newton polygon does not change when we go to an algebraic extension.

\begin{remark}
It should be emphasized that the binary operation $P_1,\:P_2\longrightarrow \mathrm{Res}_{P_1,P_2}^{V}(T)\:$ does not have all the nice properties like associativity, and that it is by no means the only one which gives a new polynomial having (almost always) $N(P_1)*N(P_2)\:$ as Newton polygon.
\end{remark}

Another feature of $N_1*N_2\:$ is that the height $h(N_1*N_2\:)$ is twice the \underline{mixed volume} of $N_1\:$ and $N_2\:$ in the following sense:

Let $N_1,\ldots,N_r$ be Newton polyhedra of finite volume in (the first quadrant of) $\mathbf{R}^d$. Then, for $\lambda_i\geq 0\:$ in $\mathbf{R}\:$ there is a polynomial expression for the volume of the polyhedron $\lambda_1N_1+\ldots+\lambda_rN_r\:$ (i.e., of the complement in the positive quadrant of the region bounded by $\sum_{i=1}^r\lambda_iN_i$):
$$\mathrm{Vol}(\lambda_1N_1+\ldots+\lambda_rN_r)=\sum_{|\alpha|=d}\frac{d!}{\alpha!}\mathrm{Vol}(N_1^{\left[\alpha_1\right]},\dots,N_r^{\left[\alpha_r\right]})\lambda_1^{\alpha_1}\ldots\lambda_r^{\alpha_r},$$\noindent
where $\alpha!=\alpha_1!\ldots\alpha_d!\:$ and $\mathrm{Vol}(N_1^{\left[\alpha_1\right]},\dots,N_r^{\left[\alpha_r\right]})\:$ is (by definition) the mixed volume of index $\alpha=(\alpha_1,\ldots,\alpha_r)\:$ of $N_1,\ldots,N_r$. We see that $\mathrm{Vol}(N_1^{\left[d\right]},N_2^{\left[0\right]},\ldots,N_r^{\left[0\right]})=\mathrm{Vol}(N_1)$, etc.

\begin{remark}
These mixed volumes should not be confused with the mixed volumes of Minkowski which are used in the isoperimetric and other external problems in the theory of convex sets. The sum of Newton polyhedra is not the Minkowski sum, but rather obtained by taking the convex hull of the Minkowski sum of the convex sets bounded by our Newton polyhedra. However, the construction is of course quite similar, and the proof of the above result is obtained by induction on $d\:$ after reduction to the case whose all the $N_i\:$ are elementary polyhedra, thanks to the following fact: let $N\:$ be a Newton polyhedron of finite volume, and let $\sigma_1,\ldots,\sigma_k\:$ be its $(d-1)$-dimensional faces. Let $h_i\:$ be the distance from the origin to the hyperplane containing $\sigma_i$. then we have:

$$d\cdot \mathrm{Vol}(N)=\sum_{i=1}^{k}h_i\mathrm{Vol}(\sigma_i),$$
where of course $\mathrm{Vol}(N)\:$ is the $d$-dimensional volume of the part of $\mathbf{R}^{d}\:$ bounded by $N\:$ in the positive quadrant, and $\mathrm{Vol}(\sigma_i)\:$ is the $(d-1)$-dimensional volume of the convex set $\sigma_i$.
\end{remark}

An interesting special case is $r=2$:
$$\mathrm{Vol}(\lambda_1N_1+\lambda_2N_2)=\sum_{i=0}^{d}{d\choose i}\mathrm{Vol}(N_1^{\left[i\right]},N_2^{\left[d-i\right]})\lambda_1^{i}\lambda_2^{d-i},$$
and in the case of polygons, $d=2$:
$$\mathrm{Vol}(\lambda_1N_1+\lambda_2N_2)=\mathrm{Vol}(N_1)\lambda_1^2+2\mathrm{Vol}(N_1^{\left[1\right]},N_2^{\left[1\right]})\lambda_1 \lambda_2+\mathrm{Vol}(N_2)\lambda_2^2.$$

\begin{proposition}\label{1.2}
Let $N_1\:$ and $N_2\:$ be two Newton polygons of finite volume. Then we have of course $\ell(N_1*N_2)=\ell(N_1) \ell(N_2)$, but also:
$$h(N_1*N_2)=2\mathrm{Vol}(N_1^{\left[1\right]}, N_2^{\left[1\right]}).$$
\end{proposition}
\noindent
Proof: Since clearly $h(N_1*N_2)\:$ is bilinear with respect to the sum of Newton polyhedra, it is sufficient to prove that $h(N_1*N_1)=2\mathrm{Vol}(N_1)$. Now if $N_1=\sum_{i=1}^{s}\Teiss{\ell_i}{h_i}\:$ it is an exercise on the convexity of the Newton polygon to check that:
$$h(N_1*N_1)=\sum_{i,j}\mathrm{min}(\ell_ih_j,\ell_jh_i)=2\mathrm{Vol}(N_1).$$

\begin{corollary}
If $N_1=\sum_{i}\Teiss{\ell_i}{h_i},\:N_2=\sum_{j}\Teiss{\ell_j'}{h_j'}\:$
$$2\mathrm{Vol}(N_1^{\left[1\right]},N_2^{\left[1\right]})=\sum_{i,j}\mathrm{min}(\ell_ih_j',\ell'_jh_i).$$
\end{corollary}

\begin{corollary}
Let $f_1,\:f_2\in\mathbf{C}\{z_1,z_2\}\:$ be such that $f_1(0,z_2),\: f_2(0,z_2)\:$ are $\neq 0$, and satisfy the condition in proposition \ref{1.2} with respect to their Newton polygons $N_1\:$ and $N_2$. Then we have:
$$\mathrm{dim}_{\mathbf{C}}(\mathbf{C}\{z_1,z_2\}/(f_1,f_2))=e(f_1,f_2)=2\mathrm{Vol}(N_1^{\left[1\right]}, N_2^{\left[1\right]}).$$
\end{corollary}
 \noindent
(This is because by proposition \ref{1.2}, and the definition of the resultant polynomial $Res^{\mathbf{C}\{z_1\}}_{f_1,f_2}(z_2)$, the height of $N_1*N_2\:$ is the valuation of the resultant, which is equal to $dim_{\mathbf{C}}(\mathbf{C}\{z_1,z_2\}/(f_1,f_2))$). (See \cite{T1}, \S 1).

\begin{remark}
It should be noted that there is no unit element in general for the product $*$. However there is an interesting subset of the set of Newton polygons, which is stable by $*\:$ and in restriction to which $*$ has a unit:
\end{remark}
\begin{definition}\label{special}
A Newton Polygon $P=\sum_{i}\Teiss{\ell_i}{h_i}\:$ is special if $\ell_i\geq h_i\:$ for all $i$.
\end{definition}\noindent
Exercise: Let $\mathbf{1}=\Teiss{1}{1}$. Show that $P\:$ is special if and only if we have $P*\mathbf{1}=P$.\par\noindent
Exercise: Enlarge the set of Newton polygons on which $*\:$ is defined to all Newton polygons, not only those of finite area, and then show that the Newton polygon one wants to write $\Teiss{1}{\infty}\:$ is a unit.\par\medskip\noindent
\centerline{\textit{The next illustration was added in 2012.}}

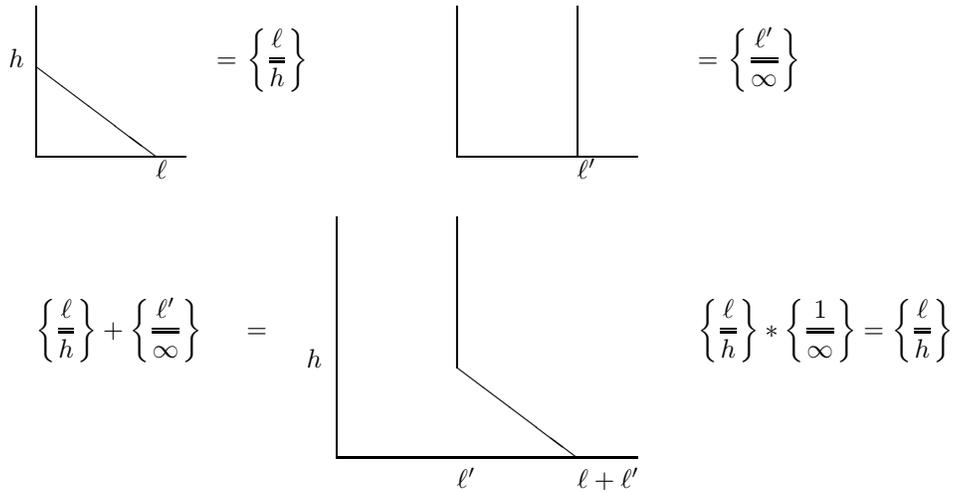
\vspace{1.5in}
\begin{picture}(1,1)\setlength{\unitlength}{.4cm}
\put(0,0){\line(0,1){5}}
\put(0,0){\line(1,0){5}}
\put(0,3){\line(4,-3){4}}
\put(-0.9,3){$h$}
\put(4,-0.7){$\ell$}
\put(6, 3){$=$}
\put(7,3){$\Teiss{\ell}{h}$}
\put(14,0){\line(0,1){5}}
\put(14,0){\line(1,0){6}}
\put(18,0){\line(0,1){5}}
\put(18,-0.7){$\ell'$}
\put(22, 3){$=$}
\put(23,3){$\Teiss{\ell'}{\infty}$}

\put(0,-6){$\Teiss{\ell}{h}+\Teiss{\ell'}{\infty}$}
\put(7,-6){$=$}
\put(9,-7){$h$}
\put(10,-10){\line(0,1){8}}
\put(10,-10){\line(1,0){10}}
\put(14, -7){\line(0,1){5}}
\put(14,-7){\line(4,-3){4}}
\put(18,-11){$\ell+\ell'$}
\put(14,-11){$\ell'$}
\put(22,-6){$\Teiss{\ell}{h}*\Teiss{1}{\infty}=\Teiss{\ell}{h}$}
\put(0,-14){ Exercise:}
\put(5,-14){$\Teiss{\ell}{h}*\Teiss{\infty}{h'}=\Teiss{\infty}{\ell h'}$,}
\put(16,-14){$\Teiss{n}{n}*\Teiss{n'}{n'}=\Teiss{nn'}{nn'}$}
\end{picture}

\vfill\eject
We mention in closing that the use of Newton polyhedra, and the use of the theory of multiplicities and integral dependence in the study of singularities are very closely related. One aspect of this is the link between the Newton polygon of two ideals, multiplicity theory and integral dependence which we saw above, but I like to mention also that there is a very close analogy between the theory of multiplicities and the theory of mixed volumes, this time in the sense of Minkowski (see \cite{E}, which suggested this to me; I thank Tadao Oda for inducing me to read it).

Here is the analogy:

\begin{center}\small
\begin{tabular}{|p{58mm}|p{64mm}|}\hline
&\\
Some Cohen-Macaulay analytic algebra $\Oc\:$ of dimension $d$&Affine space $\mathbf{R}^d$\\
&\\

&\\
An ideal $\mathbf{n}\subset\Oc$ primary for the maximal ideal& A convex body $K\subset\mathbf{R}^d\:$ of finite volume\\
&\\

&\\
Product $\mathbf{n}_1.\mathbf{n}_2$& Minkowski (pointwise) sum $K_1+K_2$\\
&\\

&\\
Multiplicity $e(\mathbf{n})$&Volume $V(K)$\\
&\\

&\\
Definition of mixed multiplicities&Definition of mixed volumes\\
$e(\mathbf{n}_1^{\left[\alpha_1\right]}+\ldots +\mathbf{n}_r^{\left[\alpha_r\right]})\:$&$V(K_1^{\left[\alpha_1\right]},\ldots,K_r^{\left[\alpha_r\right]})\:$\\
 with $\sum{\alpha_i}=|\alpha|=d\:$ by& with $\sum{\alpha_i}=|\alpha|=d\:$ by\\ 
 $e(\mathbf{n}_1^{\nu_1}\ldots \mathbf{n}_r^{\nu_r})=$ & $V(\nu_1K_1+\ldots+\nu_rK_r)=$\\
 $\sum_{|\alpha|=d}\frac{d!}{\alpha!}e(\mathbf{n}_1^{\left[\alpha_1\right]}+\ldots+\mathbf{n}_r^{\left[\alpha_r\right]})\nu_1^{\alpha_1}\ldots\nu_r^{\alpha_r}$&  $\sum_{|\alpha|=d}\frac{d!}{\alpha!}V(K_1^{\left[\alpha_1\right]},\ldots ,K_r^{\left[\alpha_r\right]})\nu_1^{\alpha_1}\ldots\nu_r^{\alpha_r}$\\
&\\

&\\
Given $\mathbf{n}_1,\:\mathbf{n}_2\:$ setting & Given $K_1,\:K_2\:$ setting\\
$e_i=e(\mathbf{n}_1^{\left[i\right]}+\mathbf{n}_2^{\left[d-i\right]})\:$ then&  $v_i=V(K_1^{\left[i\right]},K_2^{\left[d-i\right]})\:$ then\\ $\frac{e_i}{e_{i-1}}\geq \frac{e_{i-1}}{e_{i-2}}\: ,\ (2\leq i\leq d)\:$ \cite{T4}&$\frac{v_i}{v_{i-1}}\leq \frac{v_{i-1}}{v_{i-2}}\: ,\ (2\leq i\leq d)\:$ (Fenchel-Alexandrov)\\
&\\

&\\
Minkowski type inequality &Br\"unn-Minkowski inequality\\
$e(\mathbf{n}_1.\mathbf{n}_2)^{\frac{1}{d}}\leq e(\mathbf{n}_1)^{\frac{1}{d}}+e(\mathbf{n}_2)^{\frac{1}{d}}$& $V(K_1+K_2)^{\frac{1}{d}}\geq V(K_1)^{\frac{1}{d}}+V(K_2)^{\frac{1}{d}}$\\
&\\

&\\
Equality $\Longleftrightarrow \overline{\mathbf{n}_1^{a}}=\overline{\mathbf{n}_2^{b}}\: (a,\:b\in\mathbf{N})\:$ \cite{T5}&Equality $\Longleftrightarrow K_1,\:K_2\:$ are homothetic [E].\\
&\\
&\\ \hline
\end{tabular}
\end{center}
\large
\bigskip

This analogy is not too difficult to explain, in two steps:\par\noindent First take $\Oc=\mathbf{C}\{z_1,\ldots,z_d\}$. Then if an ideal $\mathbf{n}\:$ is generated by monomials, $e(\mathbf{n})=d!\mathrm{Vol}(N)\:$ where $N\:$ is the Newton Polyhedron defined by these monomials and from that it is very easy to see that if $\mathbf{n}_1,\: \mathbf{n}_2\:$ are both generated by monomials

\begin{equation*}
e(\mathbf{n}_1^{\left[i\right]},\mathbf{n}_2^{\left[d-i\right]})=d!\mathrm{Vol}(N_1^{\left[i\right]},N_2^{\left[d-i\right]}),
\end{equation*}
\noindent
where the last is the mixed volume in the sense of Newton polyhedra. Then one has to relate the mixed volumes of Newton polyhedra with the Minkowski mixed volumes of convex sets, and this is where the inequalities are reversed. The details will be given elsewhere. \par\medskip
\centerline{\textit{The references in the original list have been updated in 2012.}}

\end{document}